\documentclass[12pt]{amsart}
\usepackage{amsmath,amssymb,amsthm}
\usepackage[left=1in,right=1in,top=1.25in,bottom=1.25in]{geometry}
\usepackage[bookmarks=true,colorlinks,linkcolor=blue]{hyperref}
\usepackage{graphicx,subcaption}
\usepackage{standalone}
\usepackage{mathrsfs}
\usepackage{tikz,tikz-cd}
\usetikzlibrary{positioning,knots,patterns,hobby}
\usepackage{comment}

\providecommand{\abs}[1]{{\left|{#1}\right|}}
\providecommand{\floor}[1]{{\left\lfloor{#1}\right\rfloor}}

\newcommand{\skein}{{\mathscr{S}}}
\newcommand{\torus}{{\Sigma_{1,0}}}
\newcommand{\ptorus}{{\Sigma_{1,1}}}
\newcommand{\psphere}{{\Sigma_{0,4}}}

\newtheorem{theorem}{Theorem}[section]
\newtheorem{thm}{Theorem}
\newtheorem{lemma}[theorem]{Lemma}
\newtheorem{corollary}[theorem]{Corollary}

\newtheorem{conj}{Conjecture}

\theoremstyle{definition}

\newtheorem{remark}{Remark}[section]

\newcommand{\red}[1]{{\color{red}#1}}

\begin{document}

\title[Positivity of skein algebras]{Lower and upper bounds for positive bases of skein algebras}

\author[Thang  T. Q. L\^e]{Thang  T. Q. L\^e}
\address{School of Mathematics, 686 Cherry Street,
 Georgia Tech, Atlanta, GA 30332, USA}
\email{letu@math.gatech.edu}

\author[Dylan Thurston]{Dylan P. Thurston}
\address{Department of Mathematics, Indiana University, Bloomington, IN 47405-7106}
\email{dpthurst@indiana.edu}

\author[Tao Yu]{Tao Yu}
\address{School of Mathematics, 686 Cherry Street,
 Georgia Tech, Atlanta, GA 30332, USA}
\email{tyu70@gatech.edu}

\thanks{Supported in part by National Science Foundation. \\
2010 {\em Mathematics Classification:} Primary 57N10. Secondary 57M25.\\
{\em Key words and phrases: Kauffman bracket skein algebra, positive basis.}}

\begin{abstract}
We show that the if a sequence of normalized polynomials gives rise to a positive basis of the skein algebra of a surface, then it is sandwiched between the two types of Chebyshev polynomials. For the closed torus, we show that the normalized sequence of Chebyshev polynomials of type one $(\hat{T}_n)$ is the only one which gives a positive basis.
\end{abstract}

\maketitle

\def\Sgp{\Sigma_{g,p}}
\def\SSgp{\skein(\Sgp;R)}
\def\SSi{\skein(\Sigma;R)}
\def\BZ{{\mathbb Z}}
\def\BR{{\mathbb R}}
\section{Introduction}

\subsection{Results}
Let $R$ be a commutative integral domain with a distinguished  invertible element $q\in R$. The main examples are $R=\mathbb{Z}[q^{\pm1}]$ and $R=\mathbb{Z}$ with $q=1$. Let $\Sigma=\Sgp$ be the oriented surface of genus $g$ with $p$ points removed. The skein algebra $\SSi$ is the $R$-algebra spanned by isotopy classes of framed links in $\Sigma\times(-1,1)$ modulo the skein relation and the trivial loop relation in Figure~\ref{fig-equiv}. The product is given by superposition. For details see Section~\ref{sec.basic}.

We assume that the ring $R$ has a positive part $R_+$, see Section~\ref{sec.pos}. When $R=\mathbb{Z}$ its positive part is 
$R_+=\mathbb{Z}_+$, the set of non-negative integers, and when  $R=\mathbb{Z}[q^{\pm1}]$ its positive part is  $R_+=\mathbb{Z}_+[q^{\pm1}]$. A basis $B$ of an $R$-algebra is {\em positive} if the structure constants are in $R_+$, i.e. for any $x, y\in B$ the product $xy$ is a linear combination of elements in $B$ with coefficients in $R_+$. An important conjecture \cite{FG} in cluster algebra theory is that the skein algebra $\skein(\Sigma)$ has a positive basis. 

A sequence of polynomials $(P_n(x))_{n=0}^\infty$ with $R$ coefficients is {\em normalized} if $P_n(x)$ is monic with degree $n$ for each $n\ge0$. Note $P_0(x)=1$ by definition. A normalized sequence $(P_n)$ defines a basis $B_P$ of $\SSi$ as in Section~\ref{sec-bases}. We say a sequence of normalized polynomials $(P_n)$ is {\em positive on $\Sigma$ over $R$} if the corresponding basis $B_P$ is a positive basis of $\SSi$.

Two important sequences considered here are the normalized Chebyshev polynomials of type one $(\hat{T}_n)$ and type two $(S_n)$.
\begin{align*}
\hat{T}_0(x)&=1,&\hat{T}_1(x)&=x,&\hat{T}_2(x)&=x^2-2,&\hat{T}_n(x)&=x\hat{T}_{n-1}(x)-\hat{T}_{n-2}(x),\quad n\ge3,\\
S_0(x)&=1,&S_1(x)&=x,&&&S_n(x)&=xS_{n-1}(x)-S_{n-2}(x),\quad n\ge2.
\end{align*}

 The second named author \cite{Thurston} showed that when $R=\mathbb{Z}$ and $q=1$, the sequence $(\hat{T}_n)$ is positive on any surface,  and made a conjecture, refining an earlier one of Fock and Goncharov \cite{FG}, that $(\hat{T}_n)$ is positive when $R=\mathbb{Z}[q^{\pm1}]$ as well.
 
  The first named author \cite{LePos} showed that if the sequence $(P_n)$ is positive on a surface $\Sigma$ having genus $\ge 1$, then $(P_n) \ge (\hat{T}_n)$. Here, $(P_n)\ge (Q_n)$ means each $P_n$ is an $R_+$-linear combination of $Q_0,Q_1,\dots,Q_n$.
  
 We  extend the lower bound result in \cite{LePos} to surfaces of genus 0, and at the same time give an upper bound for all surfaces.

\begin{thm}\label{thm-main}
Suppose $R=\mathbb{Z}[q^{\pm1}]$ and $(P_n)$ is normalized with integer coefficients. 
Let $\Sigma$ be a surface with genus at least $1$ or with at least $4$ punctures. If $(P_n)$ is positive on $\Sigma$, then $(\hat{T}_n)\le(P_n)\le(S_n)$.
\end{thm}


In the case of the closed torus $\Sigma_{1,0}$, our result is much more precise.
\begin{thm}\label{thm-main2}
A normalized sequence $(P_n)$ is positive on the closed torus $\Sigma_{1,0}$ over $\BZ$ or $\BZ[q^{\pm 1}]$ if and only if $(P_n) = (\hat{T}_n)$.
\end{thm}
See Section~\ref{sec-result} for more refined statements.

This is somewhat surprising. The sequence $(S_n)$ gives simple objects in the ring of modules over the quantum group $U_q(sl_2)$. From this point of view, $(S_n)$ is more natural than the sequence $(\hat{T}_n)$. However, it does not give a positive basis on the torus. On the other hand, if the surface has negative Euler characteristic, it is conjectured  that the $(S_n)$ is positive \cite{Thurston}.
\begin{conj} Both $(\hat{T}_n)$ and $(S_n)$ are positive on any surface $\Sigma$ with negative Euler characteristic.
\end{conj} 
For related results concerning categorifications of skein modules see \cite{QW}.
 
\begin{remark}
In theorem~\ref{thm-main} we exclude the case when $(g,p)$ is one of $(0,0),(0,1), (0,2),(0,3)$ because in these cases the skein algebra $\SSgp$ is a commutative polynomial algebra and hence obviously has a positive basis; for example, the monomial basis.
\end{remark}
\begin{remark} For other types of positivity see for example \cite{MSW,LS,CKKO}.  \end{remark}

\subsection{Acknowledgments}  The authors are supported in part by NSF grants. The first author would like to thank the CIMI Excellence Laboratory,
Toulouse, France, for inviting him on a Excellence Chair during the period of January -- July 2017 when part of this work was done. 

\section{Skein algebras}\label{sec-skein}

\subsection{Basic definitions}\label{sec.basic}
\begin{figure}[h]
\centering
\begin{subfigure}[b]{0.45\linewidth}
\centering
\[
\raisebox{-0.43\height}{
\begin{tikzpicture}[scale=0.75]
\clip (0,0)circle(1);
\fill[gray!20!white] (0,0)circle(1);
\begin{knot}[clip width=5,background color=gray!20!white]
\strand (1,1)--(-1,-1);
\strand (-1,1)--(1,-1);
\end{knot}
\end{tikzpicture}
}
=
q\raisebox{-0.43\height}{
\begin{tikzpicture}[scale=0.75]
\clip (0,0)circle(1);
\fill[gray!20!white] (0,0)circle(1);
\draw (-1,-1)..controls (0,0)..(-1,1);
\draw (1,-1)..controls (0,0)..(1,1);
\end{tikzpicture}
}
+
q^{-1}\raisebox{-0.43\height}{
\begin{tikzpicture}[scale=0.75]
\clip (0,0)circle(1);
\fill[gray!20!white] (0,0)circle(1);
\draw (-1,-1)..controls (0,0)..(1,-1);
\draw (-1,1)..controls (0,0)..(1,1);
\end{tikzpicture}
}
\]
\subcaption{Skein relation}
\end{subfigure}
\hfill
\begin{subfigure}[b]{0.45\linewidth}
\centering
\[
\raisebox{-0.43\height}{
\begin{tikzpicture}[scale=0.75]
\clip (0,0)circle(1);
\fill[gray!20!white] (0,0)circle(1);
\draw (0,0)circle(.5);
\end{tikzpicture}
}
=
(-q^2-q^{-2})\raisebox{-0.43\height}{
\begin{tikzpicture}[scale=0.75]
\clip (0,0)circle(1);
\fill[gray!20!white] (0,0)circle(1);
\end{tikzpicture}
}
\]
\subcaption{Trivial loop relation}
\end{subfigure}
\caption{Defining relations for $\skein(M)$}
\label{fig-equiv}
\end{figure}
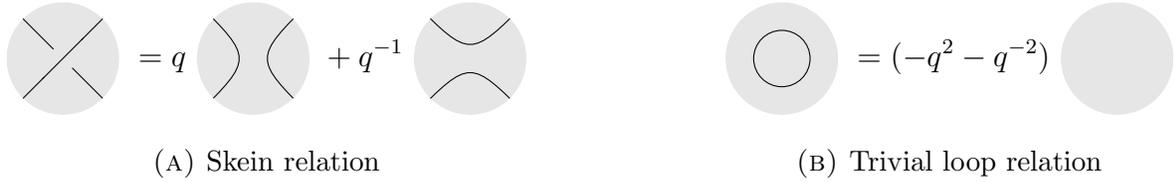

Let $R$ be a commutative integral domain with unit and with a distinguished  invertible element $q\in R$.
Suppose $M$ is an oriented $3$-manifold, not necessarily closed. A {\em framed link} $L$ in $M$ is a smooth, unoriented, closed $1$-dimensional submanifold equipped with a normal vector field. By convention, the empty set is also considered as a framed link. The skein module $\skein(M;R)$ of $M$, introduced independently by Przytycki  and Turaev \cite{Przy,Turaev1,Turaev2}, is defined as the quotient of the $R$-module freely generated by isotopy classes of framed links in $M$ by the skein relation and the trivial loop relation in Figure~\ref{fig-equiv}.

In the figure, the diagrams represent links that are identical outside of a ball in the manifold $M$. The shaded part is the projection onto the equatorial plane of the ball where the difference is. The framings on the diagrams are vertical, i.e. pointing out to the reader. The skein relations were introduced by Kauffman \cite{Kauffman}.

For $\Sigma=\Sgp$, the oriented surface with genus $g$ and $p$ punctures, let 
\[\SSi:= \skein(\Sigma\times(-1,1);R),\]
which has a product structure given by stacking. More precisely, the product of two framed links $L$ and $K$ in $\SSi$ is given by the union $i_+(L)\cup i_-(K)$, where $i_\pm:M\to M$ are the embeddings defined by $i_\pm(x,t)=(x,\frac{t\pm1}{2})$.  It is easy to see that this is a well defined product which turns  $\SSi$ into an $R$-algebra. It should be noted that there are surfaces $\Sigma\neq \Sigma'$ such that $\Sigma \times (-1,1)$ and $\Sigma' \times (-1,1)$ are diffeomorphic as 3-manifolds, but $\SSi$ and $\skein(\Sigma';R)$ have different product structures, hence different as $R$-algebras.

\newcommand{\no}[1]{ }

\no{
\begin{remark}
\red{For convenience, we replace boundary components of a surface with punctures. This does not affect the structure of the skein algebra.}
\end{remark}
}

An orientation preserving embedding of surfaces $i:\Sigma\to\Sigma'$ induces an $R$-algebra homomorphism $i_\ast:\skein(\Sigma;R)\to\skein(\Sigma';R)$ by applying $(x,t)\mapsto(i(x),t)$ to links. In particular, the mapping class group of $\Sigma$ acts on the skein algebra $\skein(\Sigma;R)$.

\def\al{\alpha}
\subsection{Multicurves}

A {\em simple multicurve} on a surface $\Sigma=\Sgp$ is a closed unoriented $1$-submanifold of $\Sigma$ none of the components of which bounds a disk in $\Sigma$. A {\em simple closed curve} is a simple multicurve with one component. A {\em peripheral curve} is simple closed curve bounding a once-punctured disk. By convention, the empty set is a simple multicurve.

A simple multicurve $\gamma$ of $\Sigma$ defines an element of $\skein(\Sigma;R)$ by the embedding $\gamma\subset\Sigma\cong\Sigma\times\{0\}$. The framing on $\gamma$ is vertical, that is, parallel to the $(-1,1)$ direction and pointing towards $1$.

If $\al$ and $\beta$ are disjoint simple multicurves, then $\al \cup \beta$ is also a simple multicurve, and $\al \beta = \beta \al = \al \cup \beta$ as elements of $\SSi$.  It follows that peripheral curves are in the center of the skein algebra, since these curves intersect trivially with any curve in the surface. Also in the skein algebra, every simple multicurve $\gamma$ can be uniquely written as
\begin{equation}\label{eqn-multi-decomp}
\gamma=\gamma_1^{n_1}\dots\gamma_r^{n_r},
\end{equation}
where $\gamma_1,\dots,\gamma_r$ are all distinct isotopy classes of the components of $\gamma$, and $n_i$ is the number of components of $\gamma$ that are isotopic to $\gamma_i$.

\def\embed{\hookrightarrow}
\subsection{Bases}\label{sec-bases}

Let $B=B(\Sigma)$ denote the set of all isotopy classes of simple multicurves. The module structure of the skein algebra is very simple.

\begin{theorem}[\cite{Przy}]\label{thm-basis}
As an $R$-module, $\SSi$ is free with  basis $B(\Sigma)$.
\end{theorem}

From this description of free bases, it is easy to see that, as an  $R$-algebra $\SSgp$ is isomorphic to a commutative polynomial algebra for the case when $g=0$ and $p \le 3$. Namely, $\skein(\Sigma_{0,0};R) \cong \skein(\Sigma_{0,1};R) \cong R$. When $(g,p)=(0,2)$, one has $\skein(\Sigma_{0,2};R)\cong R[x]$ where $x$ is the only peripheral curve. Finally $\skein(\Sigma_{0,3};R)\cong R[x,y,z]$, where $x,y,z$ are the three peripheral curves. 
 These surfaces will not be considered below. For all other $(g,p)$, over  $R=\BZ[q^{\pm 1}]$, the skein algebra $\SSgp$ is non-commutative.
 
An embedding $\iota: \Sigma \embed \Sigma'$ is {\em strict} if the induced map from $B(\Sigma)$ to $B(\Sigma')$ is injective. From Theorem~\ref{thm-basis} we get the following.
\begin{corollary}\label{r.embed}
If $\iota: \Sigma \embed \Sigma'$ is a strict embedding, then $\iota_*: \SSi \to \skein(\Sigma';R)$ is an algebra embedding.
\end{corollary}

The basis $B(\Sigma)$ can be twisted by polynomial sequences as follows. Let $P=(P_n)$ be a normalized sequence of polynomials. If $\gamma=\gamma_1^{n_1}\dots\gamma_r^{n_r}$ as in Equation \eqref{eqn-multi-decomp}, define
\[P(\gamma)=P_{n_1}(\gamma_1)\dots P_{n_r}(\gamma_r).\]
Let
\[B_P(\Sigma):=\{P(\gamma): \gamma \in B(\Sigma)\}.\]
Then $B_P(\Sigma)$ is also a free $R$-basis of $\SSi$. When $P_n(x)=x^n$ one recovers $B_P(\Sigma)=B(\Sigma)$.

\subsection{Positivity} \label{sec.pos}

Let $R_+=\mathbb{Z}_+$ if $R=\mathbb{Z}$ and $R_+=\mathbb{Z}_+[q^{\pm1}]$ if $R=\mathbb{Z}[q^{\pm1}]$. More generally, when $R$ is an arbitrary commutative domain with a distinguished invertible element $q$, a {\em positive part} of $R$ is any subset $R_+$ satisfying
\begin{enumerate}
\item $q,q^{-1}\in R_+$;
\item $R_+$ is closed under addition and multiplication;
\item $R_+\cap(-R_+)=\{0\}$.
\end{enumerate}

Fix such a positive part of $R$.
A basis of $\skein(\Sigma;R)$ is {\em positive} if the structure constants are in $R_+$, that is, for any basis elements $x,y$, the product $xy$ is an $R_+$-linear combination of the basis elements. A normalized sequence of polynomials $P=(P_n)$ is {\em positive} on $\Sigma$ over $R$ if the basis $B_P$  is positive for $\skein(\Sigma;R)$.

Recall that given two normalized sequences of polynomials  $(P_n)$ and $(Q_n)$, one defines $(P_n)\le(Q_n)$ if each $Q_n$ is an $R_+$-linear combination of $P_0,P_1,\dots,P_n$. 
\no{Since the polynomials are monic, we can always write
\[Q_n(x)=P_n(x)+\sum_{i=0}^{n-1}c_{n,i}P_i(x)\]
for some constants $c_{n,i}\in R$. We say $(P_n)\le(Q_n)$ if $c_{n,i}\in R_+$.
}

\begin{lemma} The  binary relation $(\le)$ is a partial order on the set of normalized sequences of polynomials.
\end{lemma}
\begin{proof} It is clear that $(\le)$ is reflexive and transitive, and we need to show that it is anti-symmetric. Assume that $(P_n) \le (Q_n)$ and $(Q_n) \le (P_n)$. 
Writing each sequence $(P_n)$ and $(Q_n)$ as an infinite column vectors, then $(Q_n) = A \times  (P_n)$, where $A$ is an upper triangular $\BZ_+ \times \BZ_+$ matrix, having  1 on the diagonal, and having entries $R_+$. We can write $A = I + N$, where $I$ is the identity matrix and $N$ is strictly upper triangular.  Suppose  the matrix $N= (N_{ij})$ is not 0. Among all the non-zero entries of $N$ let  $N_{ij}$ be the one with the smallest pair $(j-i,i)$ in the lexicographic order. Then it is easy to see that $(A^{-1})_{ij} = - N_{ij}$. Since 
$(Q_n) \le (P_n)$, all the entries of $A^{-1}$ are in $R_+$. It follows that both $N_{ij}$ and $-N_{ij}$ are in $R_+$, a contradiction. Thus $N=0$, and $(P_n) = (Q_n)$.
\end{proof}

The Chebyshev polynomials of type one $(T_n)$ and type two $(S_n)$ are defined by the recurrence relations
\begin{align*}
T_0(x)&=2,&T_1(x)&=x,&T_n(x)&=xT_{n-1}(x)-T_{n-2}(x),\quad n\ge2,\\
S_0(x)&=1,&S_1(x)&=x,&S_n(x)&=xS_{n-1}(x)-S_{n-2}(x),\quad n\ge2.
\end{align*}
They can be characterized by $T_n(t+t^{-1})=t^n+t^{-n}$ and $S_n(t+t^{-1})=t^n+t^{n-2}+\dots+t^{-n}$. Thus $S_n(x)-S_{n-2}(x)=T_n(x)$ for $n\ge 2$. While $(S_n)$ is a normalized sequence, $(T_n)$ is not. We normalize $T_n$ by setting $\hat{T}_0(x)=1$ and $\hat{T}_n(x)=T_n(x)$ for $n>0$, as in the introduction. Then $(\hat{T}_n)$ is a normalized sequence, and $(\hat{T}_n)\le(S_n)$ since
\[S_n(x)=\sum_{i=0}^{\floor{n/2}}\hat{T}_{n-2i}(x).\]

\def\hT{\hat T}
\subsection{Results}\label{sec-result}  Here are fuller,  more refined versions of Theorems~\ref{thm-main} and~\ref{thm-main2}. For the convenience of proofs, we formulated them in 3 statements.

\begin{theorem}\label{thm-main-lower} Suppose $R$ is a commutative domain with a distinguished invertible element $q$ and a positive part $R_+$. Assume either $g\ge 1$ or $p \ge 4$. Let $(P_n)$ be a normalized sequence of polynomials that is positive on $\Sgp$ over $R$. 
 
 (a) One has  $(P_n)\ge(\hat{T}_n)$ and $P_1(x)=x$.
 
 (b) If $P_n= \hT_n$ for $n \le 3$ then $(P_n) = (\hT_n)$.
\end{theorem}

\begin{theorem}\label{thm-main-torus} Suppose $R$ is a commutative domain with a distinguished invertible element $q$ and a positive part $R_+$. A normalized sequence polynomial  $(P_n)$  is positive on the torus $\torus$ over $R$ if and only if $(P_n)=(\hat{T}_n)$.
\end{theorem}

\begin{theorem}\label{thm-main-upper} Suppose $R=\mathbb{Z}[q^{\pm1}]$. Assume either $g\ge 1$ or $p \ge 4$. Let $(P_n)$ be a normalized sequence of polynomials integer coefficients that is positive on $\Sgp$ over $R$.  Then $(P_n)\le(S_n)$.
\end{theorem}

To prove these theorems, it suffices to consider $3$ {\em basic surfaces}: the closed torus $\Sigma_{1,0}$, the once-punctured torus $\ptorus$, and the sphere with $4$ punctures $\psphere$. 
This can be seen as follows.
If $\Sgp$ has at least $4$ punctures, then there is a strict  embedding of  $\psphere$ into $\Sgp$ and by Corollary~\ref{r.embed} the skein algebra $\skein(\psphere;R)$ embeds in $\skein(\Sigma;R)$. The results for $\psphere$ implies the corresponding results for $\Sgp$. 
If $\Sigma$ has genus at least $1$, then either $\Sigma=\torus$ or  $\ptorus$ strictly embeds into $\Sgp$. In the latter case the results for $\ptorus$  implies those for  $\Sgp$.

Theorem~\ref{thm-main-lower} is proved in Section~\ref{sec-lower}. Theorem~\ref{thm-main-torus} is proved in Section~\ref{sec-upper-torus}. Theorem~\ref{thm-main-upper} is proved in Section~\ref{sec-upper-ptorus} for $\ptorus$ and Section~\ref{sec-upper-psphere} for $\psphere$. The case of $\torus$ is a corollary of Theorem~\ref{thm-main-torus}.

\section{Parameterization of  curves on basic surfaces}\label{sec-curves}
\no{\subsection{Torus $\torus$} The free abelian group $\BZ^2$ acts on the standard plane $\BR^2$ by translation, and we identify the closed torus $\torus$ with the quotient $\BR^2/\BZ^2$. We have the canonical projection $\BR^2 \to \torus$. For each pair of coprime integers $(r,s)$ let $\tC(r,s) \subset \BR^2$ be the straight line passing through the points $(0,0)$ and $(r,s)$. Let $C(r,s) \in B$ be the {\em isotopy} class of  the projection $\tC(r,s)$. Besides let $C(0,0)$ denote the empty curve. Then the set of all $C(r,s)$, with 
}

A simple curve on the torus $\torus$ is determined by the homology class up to sign, since the curves are unoriented. After choosing a basis of homology on $\torus$, every essential simple curve can be represented by a pair of coprime integers $(r,s)$, which is identified with $(-r,-s)$. 

The isotopy classes of simple curves on $\ptorus$ except the peripheral curve are in one-to-one correspondence with those on $\torus$. Thus the same notations can be used to represent essential simple closed curves on $\ptorus$.

In both cases, the mapping class group is $SL_2(\mathbb{Z})$ and the action of a mapping class on the curves is the standard linear action.

\begin{figure}[h]
\centering
\[
\begin{tikzpicture}
\draw (0,0)ellipse[x radius=2,y radius=1];
\draw (-0.5,0)..controls (-0.2,0.2) and (0.2,0.2).. (0.5,0);
\draw (-0.8,0.2) to[curve through={(-0.5,0)..(0,-0.1)..(0.5,0)}] (0.8,0.2);
\draw (0,0)ellipse[x radius=1.2,y radius=0.6];
\draw (1.2,0)node[right]{$a$};
\draw (0,-1)arc[start angle=-90,end angle=90,x radius=0.2,y radius=0.45];
\draw[dashed] (0,-0.1)arc[start angle=90,end angle=270,
	x radius=0.2,y radius=0.45];
\draw (0,-1)node[below]{$b$};
\draw[dashed] (-3,0)--(-2.1,0) (2.1,0)--(2.6,0);
\draw[->] (-2.7,0.4)..controls (-2.5,0.2) and (-2.5,-0.2)..(-2.7,-0.4);
\draw[->] (3.5,0)--(4.5,0);
\begin{scope}[xshift=7cm]
\draw (-2,0) to[closed,curve through={(0,-1)..(2,0)..(1.2,0.4)
	..(0.5,0)..(0,-0.1)..(-0.5,0)}] (-1.2,0.4);
\fill (-2,0)circle[radius=0.05] (-0.5,0)circle[radius=0.05]
	(0.5,0)circle[radius=0.05] (2,0)circle[radius=0.05];
\draw (-1.2,0.4)arc[start angle=-180,end angle=0,
	x radius=1.2,y radius=0.8];
\draw[dashed] (-1.2,0.4)arc[start angle=-180,end angle=0,
	x radius=1.2,y radius=0.9];
\draw (1.2,0.4)node[above]{$a$};
\draw (0,-1)arc[start angle=-90,end angle=90,x radius=0.2,y radius=0.45];
\draw[dashed] (0,-0.1)arc[start angle=90,end angle=270,
	x radius=0.2,y radius=0.45];
\draw (0,-1)node[below]{$b$};
\end{scope}
\end{tikzpicture}
\]
\caption{The involution of the torus}
\label{fig-torus-inv}
\end{figure}
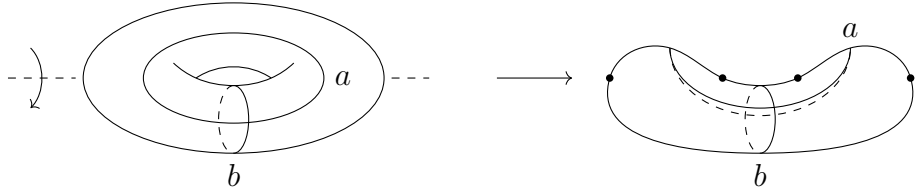

The sphere with $4$ punctures $\psphere$ can be regarded as the quotient of the torus $\torus$ by the involution in Figure~\ref{fig-torus-inv}. The action has $4$ fixed points which corresponds to the punctures. This quotient also identifies the essential simple closed curves on $\torus$ with the non-peripheral ones on $\psphere$. Thus coprime integers $(r,s)$, with the identification $(r,s)=(-r,-s)$, also represent curves on $\psphere$.

\def\al{\alpha}
\def\BZ{{\mathbb Z}}

More concretely, represent $\psphere$ as in Figure~\ref{fig-puncture-a-b}. Choose the curves $a$ and $b$, and number the punctures as in the figure. The curve surrounding puncture $i$ is denoted by $\gamma_i$. To obtain the $(r,s)$ curve on $\psphere$, take $\abs{r}$ parallel copies of $a$ and $\abs{s}$ parallel copies of $b$, and resolve each intersection such that one would turn left from $a$ to $b$ if $rs>0$ and turn right if $rs<0$. Thus $a=(1,0)$, $b=(0,1)$. The $(1,1)$ curve is demonstrated in Figure~\ref{fig-curve-b1-on-S0}. Thus the $(r,s)$ curve is the Luo product of $a^r$ and $ b^s$ \cite{Luo} if $r\ge 0$ and $s\ge 0$, and is the Luo product of $a^{-r}$ and $b^s$ if  $r<0$ and $s\ge 0$. The Luo product of two simple multicurves $\al$ and $\beta$ can be defined as the unique simple multicurve $\gamma$ such that $\al \beta = q^{I(\al,\beta)} \gamma + x$, where $I(\al,\beta)$ is the geometric intersection index of $\al$ and $\beta$, and $x$ is a $\BZ[q,q^{-1}]$-linear combination of simple multicurves with coefficients being Laurent polynomials in $q$ of highest degrees $< I(\al, \beta)$.

\begin{figure}[h]
\centering
\begin{subfigure}{0.4\linewidth}
\centering
\[
\begin{tikzpicture}
\draw[fill=gray!20!white,dotted] (0,0)ellipse[x radius=2,y radius=1];
\fill[black] (-1,0)circle[radius=0.05]
	(0,0)circle[radius=0.05]
	(1,0)circle[radius=0.05];
\draw (-1,0)node[below]{1} (0,0)node[below]{2} (1,0)node[below]{3};
\draw (-2,-0.5)node{4};
\draw (-0.5,0)ellipse[x radius=0.9,y radius=0.6];
\draw (0.5,0)ellipse[x radius=0.9,y radius=0.6];
\draw (-1.4,0)node[left]{$a$} (1.4,0)node[right]{$b$};
\end{tikzpicture}
\]
\subcaption{$a$, $b$, and puncture numbering}
\label{fig-puncture-a-b}
\end{subfigure}
\quad
\begin{subfigure}{0.4\linewidth}
\centering
\[
\begin{tikzpicture}
\draw[fill=gray!20!white,dotted] (0,0)ellipse[x radius=2,y radius=1];
\fill[black] (-1,0)circle[radius=0.05]
	(0,0)circle[radius=0.05]
	(1,0)circle[radius=0.05];
\draw (0.86,-0.14) to[closed,
	curve through={(1.2,-0.2)..(1.14,0.14)..(0.5,0.7)
	..(0,0.8)..(-0.5,0.7)..(-1.14,0.14)..(-1.2,-0.2)
	..(-0.86,-0.14)}] (0,0.4);
\end{tikzpicture}
\]
\subcaption{The $(1,1)$ curve}
\label{fig-curve-b1-on-S0}
\end{subfigure}
\caption{Curves on $\psphere$}
\label{fig-curves-on-S04}
\end{figure}

For $\psphere$, the mapping class group is $(\mathbb{Z}/2)^2\rtimes PSL_2(\mathbb{Z})$. The mapping classes that fix puncture $4$ forms a subgroup isomorphic to $PSL_2(\mathbb{Z})=SL_2(\mathbb{Z})/\{\pm1\}$. The action of this subgroup on $(r,s)$ curves is the projective linear action.

\section{Lower Bound, Proof of Theorem~\ref{thm-main-lower}}\label{sec-lower}

Note that Theorem~\ref{thm-main-lower} about lower bounds does not assume $R=\mathbb{Z}[q^{\pm1}]$.

\no{
\begin{theorem}\label{thm-lower-bound}
Let $\Sigma$ be a surface with genus at least $1$ or with at least $4$ punctures. If $(P_n)$ is positive on $\Sigma$, then $(P_n)\ge(\hat{T}_n)$. In addition, $P_1(x)=x$.
\end{theorem}}

\begin{proof}[Proof of Theorem~\ref{thm-main-lower}]

(a) 
The case when $\Sigma$ has genus at least $1$ is already proved in \cite{LePos}. Now assume $\Sigma=\psphere$. Let $\sigma$ be the counterclockwise half twist along $a$, fixing punctures $3$ and $4$ and exchanging punctures $1$ and $2$. Then $\sigma$ can be represented by the matrix $\begin{pmatrix}1&1\\0&1\end{pmatrix}$. For convenience, define $b_n=(n,1)=\sigma^n(b)$.

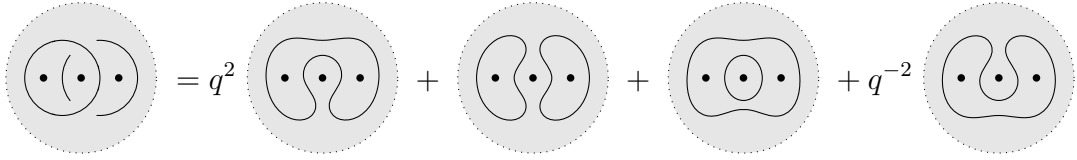
\begin{figure}[h]
\centering
\[
\raisebox{-0.43\height}{
\begin{tikzpicture}
\draw[fill=gray!20!white,dotted] (0,0)circle[radius=1];
\fill[black] (-0.5,0)circle[radius=0.05]
	(0,0)circle[radius=0.05]
	(0.5,0)circle[radius=0.05];
\begin{knot}[clip width=5,background color=gray!20!white,end tolerance=2pt]
\strand (-0.25,0)circle[radius=0.5];
\strand (0.25,0)circle[radius=0.5];
\end{knot}
\end{tikzpicture}
}
=q^2
\raisebox{-0.43\height}{
\begin{tikzpicture}
\draw[fill=gray!20!white,dotted] (0,0)circle[radius=1];
\fill[black] (-0.5,0)circle[radius=0.05]
	(0,0)circle[radius=0.05]
	(0.5,0)circle[radius=0.05];
\draw (-0.75,0) to[closed,curve through={(-0.5,-0.5)..(-0.1,-0.43)
	..(-0.25,0)..(0,0.3)..(0.25,0)..(0.1,-0.43)..(0.5,-0.5)
	..(0.75,0)..(0.5,0.5)..(0,0.5)}] (-0.5,0.5);
\end{tikzpicture}
}
+
\raisebox{-0.43\height}{
\begin{tikzpicture}
\draw[fill=gray!20!white,dotted] (0,0)circle[radius=1];
\fill[black] (-0.5,0)circle[radius=0.05]
	(0,0)circle[radius=0.05]
	(0.5,0)circle[radius=0.05];
\draw (-0.75,0) to[closed,curve through={(-0.5,-0.5)..(-0.1,-0.43)
	..(-0.25,0)..(-0.1,0.43)}] (-0.5,0.5);
\draw (0.75,0) to[closed,curve through={(0.5,-0.5)..(0.1,-0.43)
	..(0.25,0)..(0.1,0.43)}] (0.5,0.5);
\end{tikzpicture}
}
+
\raisebox{-0.43\height}{
\begin{tikzpicture}
\draw[fill=gray!20!white,dotted] (0,0)circle[radius=1];
\fill[black] (-0.5,0)circle[radius=0.05]
	(0,0)circle[radius=0.05]
	(0.5,0)circle[radius=0.05];
\draw (-0.75,0) to[closed,curve through={(-0.5,-0.5)..(-0.1,-0.43)
	..(0.1,-0.43)..(0.5,-0.5)..(0.75,0)..(0.5,0.5)..(0,0.5)}] (-0.5,0.5);
\draw (-0.25,0) to[closed,curve through={(0,0.3)..(0.25,0)}] (0,-0.3);
\end{tikzpicture}
}
+q^{-2}
\raisebox{-0.43\height}{
\begin{tikzpicture}
\draw[fill=gray!20!white,dotted] (0,0)circle[radius=1];
\fill[black] (-0.5,0)circle[radius=0.05]
	(0,0)circle[radius=0.05]
	(0.5,0)circle[radius=0.05];
\draw (-0.75,0) to[closed,curve through={(-0.5,0.5)..(-0.1,0.43)
	..(-0.25,0)..(0,-0.3)..(0.25,0)..(0.1,0.43)..(0.5,0.5)
	..(0.75,0)..(0.5,-0.5)..(0,-0.5)}] (-0.5,-0.5);
\end{tikzpicture}
}
\]
\caption{Four resolutions of $ab$}\label{fig-ab-resolve}
\end{figure}

By resolving both crossings between $a$ and $b$ as in Figure~\ref{fig-ab-resolve}, we get
\[ab=q^2b_1+q^{-2}b_{-1}+c_0,\]
where
\[c_0=\gamma_1\gamma_3+\gamma_2\gamma_4.\]
Applying $\sigma^n$ to both sides, we have
\begin{equation}\label{eqn-s04-abn}
ab_n=q^2b_{n+1}+q^{-2}b_{n-1}+c_n,
\end{equation}
where we defined
\[c_n=\sigma^n(c_0)=\begin{cases}
c_0=\gamma_1\gamma_3+\gamma_2\gamma_4,&n\text{ even},\\
c_1=\gamma_1\gamma_4+\gamma_2\gamma_3,&n\text{ odd}.
\end{cases}\]

The following can be proved easily using induction.
\begin{lemma}
For $n\ge0$,
\[T_n(a)b=q^{2n}b_n+q^{-2n}b_{-n}+c_0f_n+c_1g_n,\]
where
\[f_n(a)=\sum_{\substack{0<i\le n\\i\text{ odd}}}[i]\hat{T}_{n-i}(a)\quad\text{and}\quad g_n(a)=\sum_{\substack{0<i\le n\\i\text{ even}}}[i]\hat{T}_{n-i}(a)\]
are polynomials of $a$ with $R$ coefficients. Here $[i]=q^{2i-2}+q^{2i-6}+\dots+q^{2-2i}$ is the quantum integer.
\end{lemma}

It is also possible to express $f_n$ and $g_n$ as positive combinations of $(S_n)$. The expressions of $f_n$ and $g_n$ are not needed in the following.

First we show $P_1(x)=x$. Write $P_1(x)=x+\delta$. Then
\begin{align*}
P_1(a)P_1(b)&=(a+\delta)(b+\delta)\\
&=(q^2b_1+q^{-2}b_{-1}+c_0)+\delta(a+b)+\delta^2\\
&=q^2(P_1(b_1)-\delta)+q^{-2}(P_1(b_{-1})-\delta)+(P_1(\gamma_1)-\delta)(P_1(\gamma_3)-\delta)+\\
&\quad+(P_1(\gamma_2)-\delta)(P_1(\gamma_4)-\delta)+\delta(P_1(a)+P_1(b)-2\delta)+\delta^2.
\end{align*}
The positivity of $(P_n)$ implies that $-\delta$ and $\delta$, the coefficients of $P_1(\gamma_1)$ and $P_1(a)$ respectively, are both in $R_+$. Thus $\delta=0$, that is $P_1(x)=x$.

Now consider $P_n(x)=T_n(x)+\delta_{n-1}T_{n-1}(x)+\dots+\delta_1T_1(x)+\delta_0$ with $n\ge2$. Then
\begin{align*}
P_n(a)P_1(b)&=(T_n(a)+\delta_{n-1}T_{n-1}(a)+\dots+\delta_0)b\\
&=q^{2n}b_n+q^{-2n}b_{-n}+\delta_{n-1}(q^{2n-2}b_{n-1}+q^{2-2n}b_{1-n})+\\
&\quad+\dots+\delta_1(q^2b_1+q^{-2}b_{-1})+\delta_0b+c_0F+c_1G,
\end{align*}
where $F$ and $G$ are polynomials of $a$ with $R$ coefficients. Then by the positivity of $(P_n)$, the coefficients of $P_1(b_i)=b_i$ are in $R_+$, which implies $\delta_i\in R_+$. Thus $(P_n)\ge(\hat{T}_n)$.

(b)
Choose an essential simple closed curve $z$ on $\Sigma$ and a regular neighborhood $N$ of $z$, which is an annulus. Then $\skein(N)\cong R[z]$ is a subalgebra of $\skein(\Sigma)$, and the positivity of $(P_n)$ on $\Sigma$ implies that $(P_n(z))$ is a positive basis for $R[z]$. 

Assume $P_i(x)=\hat{T}_i(x)$ for $i<k$ where $k\ge 4$. Since $(P_n)\ge(\hat{T}_n)$, we can write
\[P_k(x)=T_k(x)+\sum_{i=0}^{k-1}\delta_i\hat{T}_i(x).\]
for some $\delta_i\in R_+$. Consider
\[P_1(z)P_{k-1}(z)=T_1(z)T_{k-1}(z)=T_k(z)+T_{k-2}(z).\]
This should be an $R_+$-linear combination of $P_0(z),\dots,P_k(z)$. The coefficient of $P_k(z)$ is $1$ by the monic condition. Thus
\[P_1(z)P_{k-1}(z)-P_k(z)=-\delta_{k-1}\hat{T}_{k-1}(z)+(1-\delta_{k-2})\hat{T}_{k-2}(z)+\sum_{i=0}^{k-3}(-\delta_i)\hat{T}_i(z)\]
is an $R_+$-linear combination of $P_0(z)=\hat{T}_0(z),\dots,P_{k-1}(z)=\hat{T}_{k-1}(z)$. This shows $\delta_i=0$ for $i<k$ except $i=k-2$. A similar argument with $P_2(z)P_{k-2}(z)$ shows $\delta_{k-2}=0$ as well. Thus $P_k(x)=\hat{T}_k(x)$. By induction, $(P_n)=(\hat{T}_n)$.
\end{proof}

\section{Upper Bound}\label{sec-upper}

The strategy to obtain an upper bound on positive polynomials is to compute the product of simple closed curves with more and more intersections.

\subsection{Proof of Theorem~\ref{thm-main-torus}}\label{sec-upper-torus}
\begin{proof}
For any pair of integers $(r,s)\ne(0,0)$, define $(r,s)_T=T_d((r/d,s/d))$ where $d=\gcd(r,s)$. For convenience, let $(0,0)_T=2$. In this notation, the basis $B_{\hat{T}}$ is
\[B_{\hat{T}}=\{(r,s)_T:(r,s)\ne(0,0)\}\cup\{1\}.\]
Note $(r,s)_T=(-r,-s)_T$. The structure constants of the skein algebra of the torus in the basis $B_{\hat{T}}$ were computed by Frohman and Gelca \cite{FrG},
\[(r,s)_T(u,v)_T=q^{rv-us}(r+u,s+v)_T+q^{-(rv-us)}(r-u,s-v)_T,\]
which shows that $(\hT_n)$ is positive.

By Theorem~\ref{thm-main-lower}, one has $(P_n)\ge(\hat{T}_n)$ and $P_1(x)=x=T_1(x)$. To show the opposite inequality, consider
\[P_1((n,1))P_1((0,1))=(n,1)_T(0,1)_T=q^n(n,2)_T+q^{-n}(n,0)_T.\]

First let $n=2$. Write $T_2(x)=P_2(x)+\delta_1P_1(x)+\delta_0$. Then
\begin{align*}
P_1((2,1))P_1((0,1))&=q^2(2,2)_T+q^{-2}(2,0)_T\\
&=q^2\left(P_2((1,1))+\delta_1P_1((1,1))+\delta_0\right)\\
&\quad+q^{-2}\left(P_2((1,0))+\delta_1P_1((1,0))+\delta_0\right).
\end{align*}
By the positivity of $(P_n(x))$, $\delta_1$ and $(q^2+q^{-2})\delta_0$ are in $R_+$. On the other hand, $P_2(x)=T_2(x)-\delta_1T_1(x)-\delta_0$ implies $-\delta_1,-\delta_0\in R_+$. Thus $\delta_1=0=\delta_0$ and $P_2(x)=T_2(x)$.

For $n>2$, $(n,2)_T$ is either $P_1((n,2))$ or $P_2((n/2,1))$. Then by the positivity of $(P_n)$, $(n,0)_T=T_n((1,0))$ is an $R_+$-linear combination of $\{P_k((1,0))\}$. Thus $(P_n)\le(\hat{T}_n)$.
\end{proof}

\subsection{Theorem~\ref{thm-main-upper}, punctured torus case}
\label{sec-upper-ptorus}
Assume $\Sigma= \Sigma_{1,1}$.
Define $T_{r,s}=T_d((r/d,s/d))$, $S_{r,s}=S_d((r/d,s/d))$ for $(r,s)\ne(0,0)$, where $d=\gcd(r,s)$. Let $T_{0,0}=S_{0,0}=1$. The products on $\ptorus$ are more complicated than the $\torus$ case. No general formula is available yet. However, to prove the theorem, only special cases are needed.

\begin{lemma}\label{lemma-pt-once}
If the curves $(r/d,s/d)$ and $(u,v)$ intersect once, where $d=\gcd(r,s)$, then
\[T_{r,s}T_{u,v}=q^{rv-su}T_{r+u,s+v}+q^{-(rv-su)}T_{r-u,s-v}.\]
\end{lemma}

\begin{proof}
By applying a mapping class, the equation can be reduced to
\[T_{d,0}T_{0,1}=q^dT_{d,1}+q^{-d}T_{d,-1}.\]
This is essentially a reformulation of Proposition 3.1 in \cite{LePos}.
\end{proof}

\begin{lemma}
On $\ptorus$, let $U$ be the peripheral curve. Then
\[T_{1,0}T_{n,2}=q^2T_{n+1,2}+q^{-2}T_{n-1,2}+(U+q^2+q^{-2})A_n,\]
where $A_n=0$ if $n$ is even, and $A_n=1$ is $n$ is odd.
\end{lemma}

\begin{proof}
If $n$ is even, then $(n/2,1)$ and $(1,0)$ intersect once. In this case, the formula is a specialization of Lemma~\ref{lemma-pt-once}. When $n=1$, the curves $(1,0)$ and $(1,2)$ intersect twice as in Figure~\ref{fig-T10T12}.

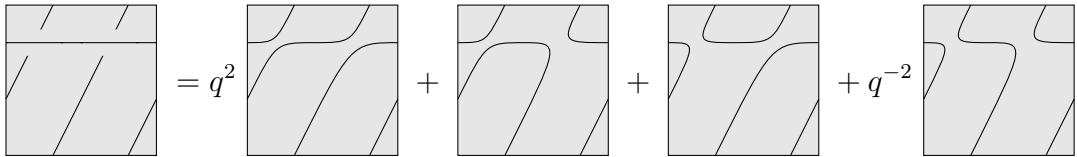
\begin{figure}[h]
\centering
\[
\raisebox{-0.43\height}{
\begin{tikzpicture}[scale=0.5]
\draw[fill=gray!20!white] (0,0)--(0,4)--(4,4)--(4,0)--cycle;
\begin{knot}[clip width=5,background color=gray!20!white,end tolerance=2pt]
\strand (0,3)--(4,3);
\strand (0,1.5)--(1.25,4) (1.25,0)--(3.25,4) (3.25,0)--(4,1.5);
\end{knot}
\end{tikzpicture}
}
=q^2
\raisebox{-0.43\height}{
\begin{tikzpicture}[scale=0.5]
\draw[fill=gray!20!white] (0,0)--(0,4)--(4,4)--(4,0)--cycle;
\draw (0,3)..controls(0.75,3)..(1.25,4);
\draw (3.25,4)..controls(2.75,3)..(1.75,3)..controls(0.75,3)..(0,1.5);
\draw (1.25,0)..controls(2.75,3)..(4,3);
\draw (3.25,0)--(4,1.5);
\end{tikzpicture}
}
+
\raisebox{-0.43\height}{
\begin{tikzpicture}[scale=0.5]
\draw[fill=gray!20!white] (0,0)--(0,4)--(4,4)--(4,0)--cycle;
\draw (0,3)..controls(0.75,3)..(1.25,4);
\draw (1.25,0)..controls(2.75,3)..(1.75,3)..controls(0.75,3)..(0,1.5);
\draw (3.25,4)..controls(2.75,3)..(4,3);
\draw (3.25,0)--(4,1.5);
\end{tikzpicture}
}
+
\raisebox{-0.43\height}{
\begin{tikzpicture}[scale=0.5]
\draw[fill=gray!20!white] (0,0)--(0,4)--(4,4)--(4,0)--cycle;
\draw (0,1.5)..controls(0.75,3)..(0,3);
\draw (3.25,4)..controls(2.75,3)..(1.75,3)..controls(0.75,3)..(1.25,4);
\draw (1.25,0)..controls(2.75,3)..(4,3);
\draw (3.25,0)--(4,1.5);
\end{tikzpicture}
}
+q^{-2}
\raisebox{-0.43\height}{
\begin{tikzpicture}[scale=0.5]
\draw[fill=gray!20!white] (0,0)--(0,4)--(4,4)--(4,0)--cycle;
\draw (0,1.5)..controls(0.75,3)..(0,3);
\draw (1.25,0)..controls(2.75,3)..(1.75,3)..controls(0.75,3)..(1.25,4);
\draw (3.25,4)..controls(2.75,3)..(4,3);
\draw (3.25,0)--(4,1.5);
\end{tikzpicture}
}
\]
\caption{Four resolutions of $T_{1,0}T_{1,2}$}\label{fig-T10T12}
\end{figure}

Opposite sides of the squares are identified, and the corners represent the puncture. Resolving both crossings, we get
\begin{align*}
T_{1,0}T_{1,2}&=q^2(1,1)^2+U+(-q^2-q^{-2})+q^{-2}(0,1)^2\\
&=q^2T_{2,2}+q^{-2}T_{0,2}+U+q^2+q^{-2}.
\end{align*}

For a general odd $n$, apply $(n-1)/2$ Dehn twists along $(1,0)$ to the equation above gives
\[T_{1,0}T_{n,2}=q^2T_{n+1,2}+q^{-2}T_{n-1,2}+(U+q^2+q^{-2}).\qedhere\]
\end{proof}

\begin{lemma}
On $\ptorus$,
\[T_{n,1}T_{0,1}=q^nT_{n,2}+q^{-n}T_{n,0}+(U+q^2+q^{-2})G_n,\]
where $G_0=G_1=0$, and
\[G_n
=\sum_{i=1}^{\floor{n/2}}q^{4i-n-2}S_{n-2i,0}\]
if $n>1$.
\end{lemma}

\begin{proof}
When $n=0,1$, the result follows from direct calculation. To use induction, note $(n,1)$ and $(1,0)$ intersect at one point. Thus $T_{1,0}T_{n,1}=qT_{n+1,1}+q^{-1}T_{n-1,1}$.
\begin{align*}
T_{n+1,1}T_{0,1}&=(q^{-1}T_{n,1}T_{1,0}-q^{-2}T_{n-1,1})T_{0,1}\\
&=q^{-1}\left(q^nT_{n,2}+q^{-n}T_{n,0}+(U+q^2+q^{-2})G_n\right)T_{1,0}\\
&\quad-q^{-2}\left(q^{n-1}T_{n-1,2}+q^{1-n}T_{n-1,0}+(U+q^2+q^{-2})G_{n-1}\right)\\
&=q^{n-1}(q^2T_{n+1,2}+q^{-2}T_{n-1,2}+(U+q^2+q^{-2})A_n)\\
&\quad+q^{-n-1}(T_{n+1,0}+T_{n-1,0})+(U+q^2+q^{-2})q^{-1}G_nT_{1,0}\\
&\quad-q^{n-3}T_{n-1,2}-q^{-n-1}T_{n-1,0}-(U+q^2+q^{-2})q^{-2}G_{n-1}\\
&=q^{n+1}T_{n+1,2}+q^{-n-1}T_{n+1,0}+(U+q^2+q^{-2})(q^{-1}G_nT_{1,0}-q^{-2}G_{n-1}+q^{n-1}A_n)\\
&=q^{n+1}T_{n+1,2}+q^{-n-1}T_{n+1,0}+(U+q^2+q^{-2})G_{n+1},
\end{align*}
where the last equality can be directly verified using the expression of $G_n$. Thus the formula holds by induction.
\end{proof}

\begin{proof}[Proof of Theorem~\ref{thm-main-upper}, punctured torus case]
Since $(P_n)$ is positive, $P_1(t)=t$. Thus
\[P((n,1))P((0,1))=T_{n,1}T_{0,1}\]
is an $R_+$-linear combination of basis elements in $B_P$. Since $(P_n)$ has integer coefficients, terms with different exponents of $q$ are separately $\mathbb{Z}_+$-linear combinations of basis elements in $B_P$. Rearranging the product in terms of the exponents of $q$, we have
\[P((n,1))P((0,1))=q^{-n}S_n((1,0))+(\text{higher degree in }q).\]
Therefore, $(P_n)\le(S_n)$.
\end{proof}

\subsection{Theorem~\ref{thm-main-upper}, sphere with $4$ punctures case}
\label{sec-upper-psphere}
The proof is similar to the punctured torus case. Define $S_{r,s}=S_d((r/d,s/d))$ for $(r,s)\ne(0,0)$, where $d=\gcd(r,s)$. Let $S_{0,0}=1$.

Recall that $\gamma_i$ is the curve surrounding puncture $i$, and $c_0=\gamma_1\gamma_3+\gamma_2\gamma_4$. $\sigma=\begin{pmatrix}1&1\\0&1\end{pmatrix}$ represents the counterclockwise half twist along $a$. $c_k=\sigma^k(c_0)$.

\begin{lemma}\label{lemma-S10Sn2}
On $\psphere$,
\begin{align*}
S_{1,0}S_{2k,2}&=q^4S_{2k+1,2}+q^{-4}S_{2k-1,2}+c_0S_{k,1}+[S_{1,0}+(q^2+q^{-2})(\gamma_1\gamma_2+\gamma_3\gamma_4)],\\
S_{1,0}S_{2k+1,2}&=q^4S_{2k+2,2}+q^{-4}S_{2k,2}+(q^2c_kS_{k+1,1}+q^{-2}c_{k+1}S_{k,1})+\Gamma.
\end{align*}
where $\Gamma=\gamma_1\gamma_2\gamma_3\gamma_4+\gamma_1^2+\gamma_2^2+\gamma_3^2+\gamma_4^2-2$.
\end{lemma}

\begin{proof}
By direct computation,
\begin{align*}
S_{1,0}S_{0,2}&=q^4S_{1,2}+q^{-4}S_{-1,2}+c_0S_{0,1}+[S_{1,0}+(q^2+q^{-2})(\gamma_1\gamma_2+\gamma_3\gamma_4)],\\
S_{1,0}S_{1,2}&=q^4S_{2,2}+q^{-4}S_{0,2}+(q^2c_0S_{1,1}+q^{-2}c_1S_{0,1})+\Gamma.
\end{align*}
The equations follow by applying $\sigma^k$.
\end{proof}

\begin{lemma}
On $\psphere$, for $n\ge0$,
\[S_{n,1}S_{0,1}=q^{2n}S_{n,2}+q^{-2n}S_{n,0}+g_n+h_n.\]
Here, $g_0=g_1=0$, and for $n\ge 2$,
\[g_n=\sum_{i=1}^{\floor{n/2}}q^{4i-2}\sum_{j=i}^{n-i}c_{n-j+1}S_{j,1}.\]
$h_n$ is a polynomial of $q$, $S_{1,0}$ and $\gamma_i$'s only. $h_0=0$, and for $n\ge1$, the exponents of $q$ in $h_n$ are between $-2n+2$ and $2n-2$.
\end{lemma}

\begin{proof}
When $n=0$, the equation clearly holds. When $n=1$, resolving both crossings of $S_{1,1}$ and $S_{0,1}$ yields
\[S_{1,1}S_{0,1}=q^2S_{1,2}+q^{-2}S_{1,0}+(\gamma_1\gamma_2+\gamma_3\gamma_4).\]
Thus the equation holds with $h_1=\gamma_1\gamma_2+\gamma_3\gamma_4$.

When $n>1$,
\begin{align*}
S_{n+1,1}S_{0,1}&=q^{-2}S_{1,0}(S_{n,1}S_{0,1})-q^{-4}S_{n-1,1}S_{0,1}-q^{-2}c_nS_{0,1}\\
&=q^{-2}S_{1,0}(q^{2n}S_{n,2}+q^{-2n}S_{n,0}+g_n+h_n)\\
&\quad-q^{-4}(q^{2n-2}S_{n-1,2}+q^{2-2n}S_{n-1,0}+g_{n-1}+h_{n-1})-q^{-2}c_nS_{0,1}.\\
&=q^{2n-2}S_{1,0}S_{n,2}+q^{-2n-2}(S_{n+1,0}+S_{n-1,0})+q^{-2}S_{1,0}g_n+q^{-2}S_{1,0}h_n\\
&\quad-q^{2n-6}S_{n-1,2}-q^{-2n-2}S_{n-1,0}-q^{-4}g_{n-1}-q^{-4}h_{n-1}-q^{-2}c_nS_{0,1}\\
&=(q^{2n-2}S_{1,0}S_{n,2}-q^{2n-6}S_{n-1,2})+q^{-2n-2}S_{n+1,0}\\
&\quad+(q^{-2}S_{1,0}g_n-q^{-4}g_{n-1}-q^{-2}c_nS_{0,1})+(q^{-2}S_{1,0}h_n-q^{-4}h_{n-1}).
\end{align*}
To continue, apply Lemma~\ref{lemma-S10Sn2} to the first term. For the product $S_{1,0}g_n$, we can use Equation~\ref{eqn-s04-abn}, written in the notations of this section
\[S_{1,0}S_{n,1}=q^2S_{n+1,1}+q^{-2}S_{n-1,1}+c_n.\]
After a routine reduction, the product $S_{n,1}S_{0,1}$ has the desired form.
\end{proof}

\begin{proof}[Proof of Theorem~\ref{thm-main-upper}, sphere with $4$ punctures case]
In the product $(n,1)(0,1)=S_{n,1}S_{0,1}$, the terms with the lowest $q$-exponent is
\[q^{-2n}S_{n,0}=q^{-2n}S_n((1,0))\]
for $n>0$. This shows that $S_n((1,0))$ is a $\mathbb{Z}_+$-linear combination of $\{P_k((1,0))\}$. Thus $(P_n)\le(S_n)$. This completes the proof of Theorem~\ref{thm-main-upper}.
\end{proof}

\end{document}